\documentclass[11pt]{article}
\usepackage{latexsym,amsmath,color,amsthm,amssymb,epsfig,graphicx,mathrsfs}
\usepackage{graphicx}
\usepackage{amssymb}
\usepackage[left=1in,top=1in,right=1in,bottom=1in]{geometry}
\usepackage[linktocpage=true]{hyperref}
\usepackage{setspace}
\usepackage{amssymb, amsmath, amsthm, graphicx,mathrsfs}
\usepackage{caption}

\usepackage{comment}
\def\qed{\hfill\ifhmode\unskip\nobreak\fi\quad\ifmmode\Box\else\hfill$\Box$\fi}
\def\ite#1{\hfill\break${}$\hbox to 50pt {\quad(#1)\hfill}}
\def\cA{{\mathcal A}}
\def\cB{{\mathcal B}}

\def\cD{{\mathcal D}}

\def\cF{{\mathcal F}}

\def\cH{{\mathcal H}}
\def\cM{{\mathcal M}}
\def\cS{{\mathcal S}}

\def\cV{{\mathcal V}}
\def\EG{{\rm EG}}

\def\SDRP{{\rm SDRP }}
\newtheorem{thm}{Theorem}[section]

\newtheorem{const}[thm]{Construction}
\newtheorem{definition}[thm]{Definition}

\newtheorem{lem}[thm]{Lemma}

\newtheorem{claim}[thm]{Claim}

\begin{document}

\pagestyle{myheadings}
\markright{{\small{\sc Z.~F\"uredi, A.~Kostochka, and Ruth Luo:   Avoiding long Berge cycles II    
}}}

\title{\vspace{-0.5in} Avoiding long Berge cycles II, exact bounds for all $n$
}

\author{
{{Zolt\'an F\" uredi}}\thanks{
\footnotesize {Alfr\' ed R\' enyi Institute of Mathematics, Hungary.
E-mail:  \texttt{z-furedi@illinois.edu}. 
Research supported in part by the Hungarian National Research, Development and Innovation Office NKFIH grant K116769, and  by the Simons Foundation Collaboration Grant 317487.
}}
\and
{{Alexandr Kostochka}}\thanks{
\footnotesize {University of Illinois at Urbana--Champaign, Urbana, IL 61801
 and Sobolev Institute of Mathematics, Novosibirsk 630090, Russia. E-mail: \texttt {kostochk@math.uiuc.edu}.
 Research 
is supported in part by NSF grant  DMS-1600592
and grants 18-01-00353A and 16-01-00499  of the Russian Foundation for Basic Research.
}}
\and{{Ruth Luo}}\thanks{University of Illinois at Urbana--Champaign, Urbana, IL 61801, USA. E-mail: {\tt ruthluo2@illinois.edu}. Research of this author
is supported in part by Award RB17164 of the Research Board of the University of Illinois at Urbana-Champaign.
}}

\date{ \today}

\maketitle

\vspace{-0.3in}

\begin{abstract}
Let   $\EG_r(n,k)$ denote the maximum number of edges in
 an $n$-vertex $r$-uniform hypergraph with no Berge cycles of length $k$ or longer. 
 In the first part of this work~\cite{FKL}, we have found exact values of $\EG_r(n,k)$ and described the structure of extremal
hypergraphs for the case when $k-2$ divides $n-1$ and   $k\geq r+3$.

In this paper we determine $\EG_r(n,k)$ and describe the extremal hypergraphs  for all $n$ when $k\geq r+4$.

\medskip\noindent
{\bf{Mathematics Subject Classification:}} 05D05, 05C65, 05C38, 05C35.\\
{\bf{Keywords:}} Berge cycles, extremal hypergraph theory.
\end{abstract}

\section{ \bf  Definitions, Berge $F$ subhypergraphs}

An $r$-uniform hypergraph, or simply {\em $r$-graph}, is a family of $r$-element subsets of a finite set.
We associate an $r$-graph $\cH$ with its edge set and call its vertex set $V(\cH)$.
Usually we take $V(\cH)=[n]$, where $[n]$ is the set of first $n$ integers, $[n]:=\{ 1, 2, 3,\dots, n\}$.
We also use the notation $\cH\subseteq \binom{[n]}{r}$.

\begin{definition}
For a graph $F$ with vertex set $\{v_1, \ldots, v_p\}$ and edge set $\{e_1, \ldots, e_q\}$, a hypergraph $\mathcal H$ contains a {\bf Berge $F$} if there exist distinct vertices $\{w_1, \ldots, w_p\} \subseteq V(\mathcal H)$ and  edges $\{f_1, \ldots, f_q\} \subseteq E(\mathcal H)$, such that if $e_i = v_{\alpha} v_{\beta}$, then $\{w_{\alpha}, w_{\beta}\} \subseteq f_i$. 
\end{definition}

Of particular interest to us are Berge cycles and Berge paths.

\begin{definition} A {\bf Berge cycle} of length $\ell$ in a hypergraph is a set of $\ell$ distinct vertices $\{v_1, \ldots, v_\ell\}$ and $\ell$ distinct edges $\{e_1, \ldots, e_\ell\}$ such that $\{ v_{i}, v_{i+1} \}\subseteq   e_i$ with indices taken modulo $\ell$.

  A {\bf Berge path} of length $\ell$ in a hypergraph  is a set of $\ell+1$ distinct vertices $\{v_1, \ldots, v_{\ell+1}\}$ and $\ell$ distinct hyperedges $\{e_1, \ldots, e_{\ell}\}$ such that $\{ v_{i}, v_{i+1} \}\subseteq   e_i$ for all $1\leq i\leq \ell$.
\end{definition}

Let $\cH$ be a hypergraph  and  $p$ be an integer. The {\em $p$-shadow}, $\partial_p \cH$, is the collection of the $p$-sets that lie in some edge of $\cH$. In particular, we will often consider the 2-shadow $\partial_2 \cH$ of a $r$-uniform hypergraph $\cH$.  Each edge of $\cH$ yields in  $\partial_2 \cH$ a clique on $r$ vertices.

\section{Graphs without long cycles}\label{sec2}

\begin{thm}[Erd\H{o}s and Gallai~\cite{ErdGal59}]\label{EGpaths}
Let $k\geq 3$ and let $G$ be an $n$-vertex graph with no cycle of length $k$ or longer. Then $e(G) \leq (k-1)(n-1)/2$. 
\end{thm}
This bound is the best possible if $n-1$ is divisible by $k-2$.
A matching lower bound can be obtained by gluing together complete graphs of sizes $k-1$.

Let $\EG(n,k)$ denote the maximum size of a graph on $n$ vertices such that it does not contain any cycle of length $k$ or longer.
Write $n$ in the form of $(k-2)\lfloor \frac{n-1}{k-2}\rfloor +m$ where $1\leq m\leq k-2$.
Considering an $n$-vertex graph whose $2$-connected blocks are complete graphs of size $k-1$ except one which is a $K_m$ we get

\begin{equation}\label{eq:1}
  \EG(n,k) \geq  f(n,k):=  \left\lfloor \frac{n-1}{k-2} \right\rfloor \binom{k-1}{2} +   \dbinom{m}{2}.
\end{equation}

It took some 15 years to prove that  equality holds in~\eqref{eq:1} for all $n$ and $k\geq 3$ (Kopylov~\cite{Kopy} and independently Woodall~\cite{WoodallA}).
One of the difficulties is, as Faudree and Schelp~\cite{FaudScheB, FaudSche75} observed, that for odd  $k$ there are
 infinitely many extremal graphs very different from the ones above.

\begin{const}\label{const:22}
Fix $k\geq 4$, $n \geq k$, $\frac{k}{2} > a\geq 1$. Define the $n$-vertex graph $H_{n,k,a}$ as follows.
 The vertex set of $H_{n,k,a}$ is partitioned into three sets $A,B,C$ such that $|A| = a$, $|B| = n - k + a$ and $|C| = k - 2a$
 and the edge set of $H_{n,k,a}$ consists of all edges between $A$ and $B$ together with all edges in $A \cup C$.
$B$ is taken to be an independent set. 
   \end{const}

\begin{figure}[!ht]
  \centering
    \includegraphics[width=0.25\textwidth]{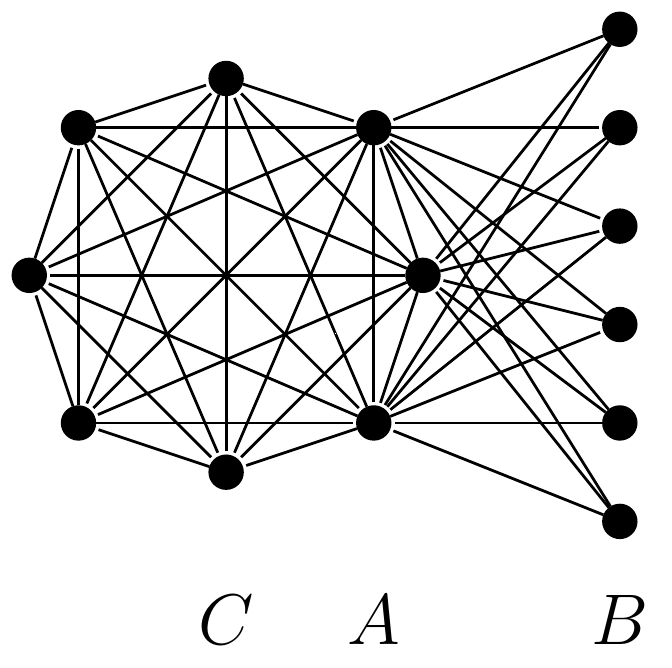}
  \caption{$H_{14,11,3}$}
\end{figure}

When $a \geq 2$, $H_{n,k,a}$ is 2-connected, has no cycle of length $k$ or longer, and
\[ e(H_{n,k,a}) =\binom{k-a}{2}+  a(n-k+a).\]

Kopylov and Woodall   (\cite{Kopy} and~\cite{WoodallA}) characterized the structure of  the extremal graphs. Namely, either \\
--- the blocks of $G$ are $p$ complete graphs $K_{k-1}$ and a $K_m$,
  where $p:=\lfloor \frac{n-1}{k-2}\rfloor$, or
\\
--- $k$ is odd, $m=(k+1)/2$ or $(k-1)/2$ and $q$ of the blocks of $G$ are
 $K_{k-1}$'s and one block is a copy of an $H_{n-q(k-2), k, (k-1)/2}$.

\section{\bf  Main result: Hypergraphs with no long Berge cycles}\label{sec3}

Let $\EG_r(n,k)$ denote the maximum size of an $r$-uniform hypergraph on $n$ vertices that does not contain any Berge cycle of length $k$ or longer.
In~\cite{FKL}, we proved an analogue of the Erd\H{o}s--Gallai theorem on cycles for $r$-graphs.

\begin{thm}[\cite{FKL}]\label{mainF}Let $r \geq 3$ and $k \geq r+3$, and suppose $\mathcal H$ is an  $n$-vertex $r$-graph with no Berge cycle of length $k$ or longer. Then $e(\mathcal H) \leq \frac{n-1}{k-2}{k-1 \choose r} $. Moreover, equality is achieved if and only if $\partial_2 \cH$ is connected and for every block $D$ of $\partial_2 \cH$, $D = K_{k-1}$ and $\mathcal H[D] = K_{k-1}^{(r)}$.
\end{thm}

\begin{center}
\includegraphics[scale=.4]{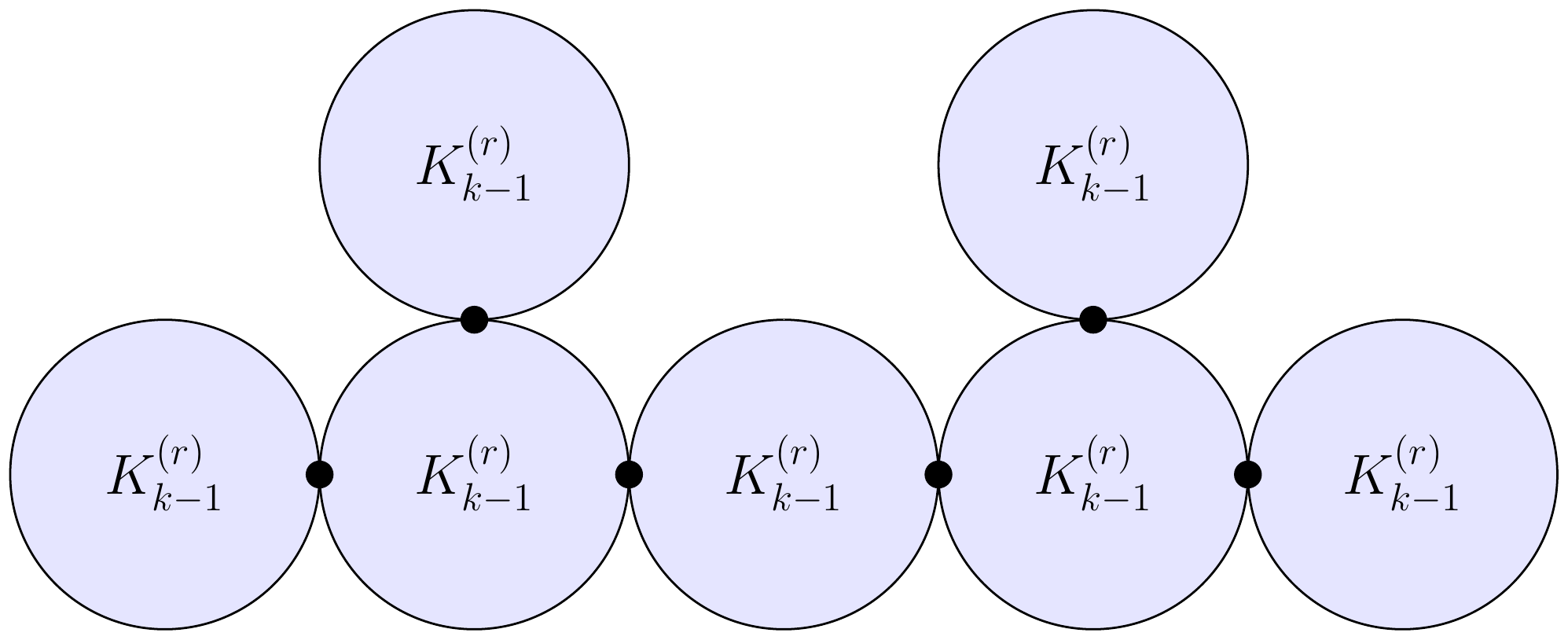}
\end{center}

Since a Berge cycle can only be contained in a single block of the 2-shadow $\partial_2\cH$,
 the construction in Theorem~\ref{mainF} cannot contain Berge cycles of length $k$ or longer.
Thus Theorem~\ref{mainF} determines $\EG_r(n,k)$ and describes extremal $r$-graphs  when $k-2$ divides $n-1$ and   $k\geq r+3$.
 Ergemlidze,  Gy\H ori,  Methuku,   Salia,  Tompkins, and  Zamora~\cite{EGMSTZ} proved  similar results
 for $k\in \{r+1,r+2\}$. 
 The case of short cycles, $k\leq r$,  is different, see~\cite{KosLuo}.

Our goal in this paper is to determine $\EG_r(n,k)$ {\em for all $n$} when $r\geq 3$ and $k \geq r+4$.
We also describe the extremal hypergraphs.
We {\em conjecture} that our results below holds for $k=r+3$ too. The tools used here do not seem to be sufficient to verify the conjecture (see the remark at the end of Section~\ref{sec:2r}).
The case $n\leq k-1$ is trivial, $\EG_r(n,k)={n \choose r}$.
Let $n=(k-2)\lfloor \frac{n-1}{k-2}\rfloor +m$ where $1\leq m\leq k-2$.
Define

\begin{equation}\label{eq:32}
   f_r(n,k):=  \left\lfloor \frac{n-1}{k-2} \right\rfloor \binom{k-1}{r} +
\left\{\begin{array}{ll}m-1         &\mbox{for  $1\leq m\leq r$},\\
   \dbinom{m}{r}   &\mbox{for $r+1\leq m \leq k-2$}. \end{array} \right.
\end{equation}

\begin{thm}\label{main:all}Let $r \geq 3$ and $k \geq r+4$, and suppose $\mathcal H$ is an  $n$-vertex $r$-graph with no Berge cycle of length $k$ or longer. Then $e(\mathcal H) \leq f_r(n,k) $. Moreover, equality is achieved if and only if $\cH$ has the structure described in Constructions~{\rm \ref{const41}} and~{\rm \ref{const42}}  in the next section.
\end{thm}

{\em The structure of the paper is as follows.}\\
In the next section (Section~\ref{sec:const}) we prove the lower bound $\EG_r(n,k)\geq f_r(n,k)$.
In Section~\ref{sec:tools} we recall some tools we developed in~\cite{FKL}: the notion of {\em representative pairs} and Kopylov's Theorem in a useful form.
In Section~\ref{sec:2r} we introduce one more tool, the notion of $(2,r)$ mixed hypergraphs and propose a more general problem.
In Section~\ref{sec:inequalities} we prepare the proof by proving a handy upper bound in the case of 
a $2$-connected $\partial_2\cH$, and finally
 in Section~\ref{sec:proof} we prove our main result, Theorem~\ref{main:all}.

\section{Constructions} \label{sec:const}

In this section we define two classes of $r$-graphs avoiding Berge cycles of length $k$ or longer (for $k\geq r+2$).
Write $n$ in the form of $(k-2)\lfloor \frac{n-1}{k-2}\rfloor +m$ where $1\leq m\leq k-2$.
Let $p:=\lfloor \frac{n-1}{k-2}\rfloor$.
Let $V=[n]$ be an $n$-element set (the set of vertices).

\begin{const}\label{const41}
In case of $m\geq r+1$, let $V_1, \dots, V_{p+1}$ be a sequence of subsets of $[n]$ satisfying
\begin{equation}\label{eq:const41}
    | (V_1 \cup \dots \cup V_{i-1}) \cap V_{i}|= 1,
\end{equation}
 for all $1 < i \leq p+1$ such that one $V_i$ has $m$ elements and  each other  has $(k-1)$-elements. 
Then replace each $V_i$ with a copy of $K^{(r)}_{|V_i|}$, the complete $r$-uniform hypergraph on it.
    \end{const}

Each Berge cycle in the $r$-uniform families in Construction~\ref{const41} must be contained in one of the $V_i$'s so its length is at most $k-1$. Hence
\[ 
   \EG_r(n,k)\geq     p \binom{k-1}{r} + \dbinom{m}{r}
\]
for all $n$, $k$, and $r$. We will see in Section~\ref{sec:proof} that in case of $m\geq r+1$  (and $k\geq r+4\geq 7$) these are the only extremal hypergraphs.

\begin{const}\label{const42}
In case of $m\leq r$, let $\cV:= \{ V_1, \dots, V_p\}$  be a sequence of $(k-1)$-element subsets of $[n]$ such that
\begin{equation}\label{eq:const42}
    | (V_1 \cup \dots \cup V_{i-1}) \cap V_{i}|\leq 1
\end{equation}
for every $i\geq 2$.
Let $H$ be the graph whose edge set is the union of the edge sets of complete graphs on $V_i\in \cV$, so $|E(H)|=p\binom{k-1}{2}$.
Then $H$ has a forest-like structure of cliques
 (i.e., every block of $H$ is a clique), and in particular every cycle is contained in some $V_i \in \cV$.

The graph $H$ necessarily consists of $m$ (nonempty) components, let $C_1, \dots, C_m$ be the vertex sets of them. 
Some $C_\alpha$'s could be singletons, and $\bigcup_{\alpha=1}^m C_\alpha=V$. 
Let $H_\alpha:= H|C_\alpha$.
Define $\cB_i$ as the complete $r$-graph with vertex set $V_i$, and set
$\cH_\alpha:= \cup \{ \cB_i: V_i\in \cV, V_i\subset C_\alpha\}$, $\cH := \cup_{\alpha=1}^m \cH_\alpha$. 
%

If $m > 1$, let $T$ be a tree with vertex set $[m]$ such that  a pair $e=\{\alpha(e), \alpha'(e)\}$ is in $E(T)$ only if the components $C_\alpha$ and $C_{\alpha'}$ in $H$ satisfy $|V(C_\alpha)| + |V(C_{\alpha'})| \geq r$. For each such edge $e$, we ``blow up" $e$ into an $r$-edge containing vertices of $C_\alpha$ and $C_{\alpha'}$ as follows:

Select the non-empty sets $A(e)\subseteq C_\alpha$ and  $A'(e)\subseteq C_{\alpha'}$ so that $|A(e)|+ |A'(e)|=r$ and if $|V(C_\alpha)| > 1$ (resp. $|V(C_{\alpha'})| > 1$), then   $A(e)\subseteq V_i \subseteq C_\alpha$ for some $V_i\in \cV$ (resp. $A'(e)\subseteq V_{i'} \subseteq C_{\alpha'}$ for some $V_{i'}\in \cV$). 
Let $\cD:=\{ A(e) \cup A'(e): e\in E(T) \}$. Our construction is $\cH\cup \cD$ (see Figure 2).

 \end{const}

By definition, $\cH\cup \cD$   has no long Berge cycle yielding
\[
   \EG_r(n,k)\geq |\cH|+|\cD|  =      p \binom{k-1}{r} + m-1
   \]
for all $n$, $k$, and $r$.
Indeed, every edge of $\cD$ is a {\em cut-edge} of the hypergraph $\cH\cup \cD$, every Berge cycle of  $\cH\cup \cD$ is contained in a single component $C_\alpha$, even more, it is contained a single $V_i$.

\begin{figure}
\centering
\includegraphics[scale=.4]{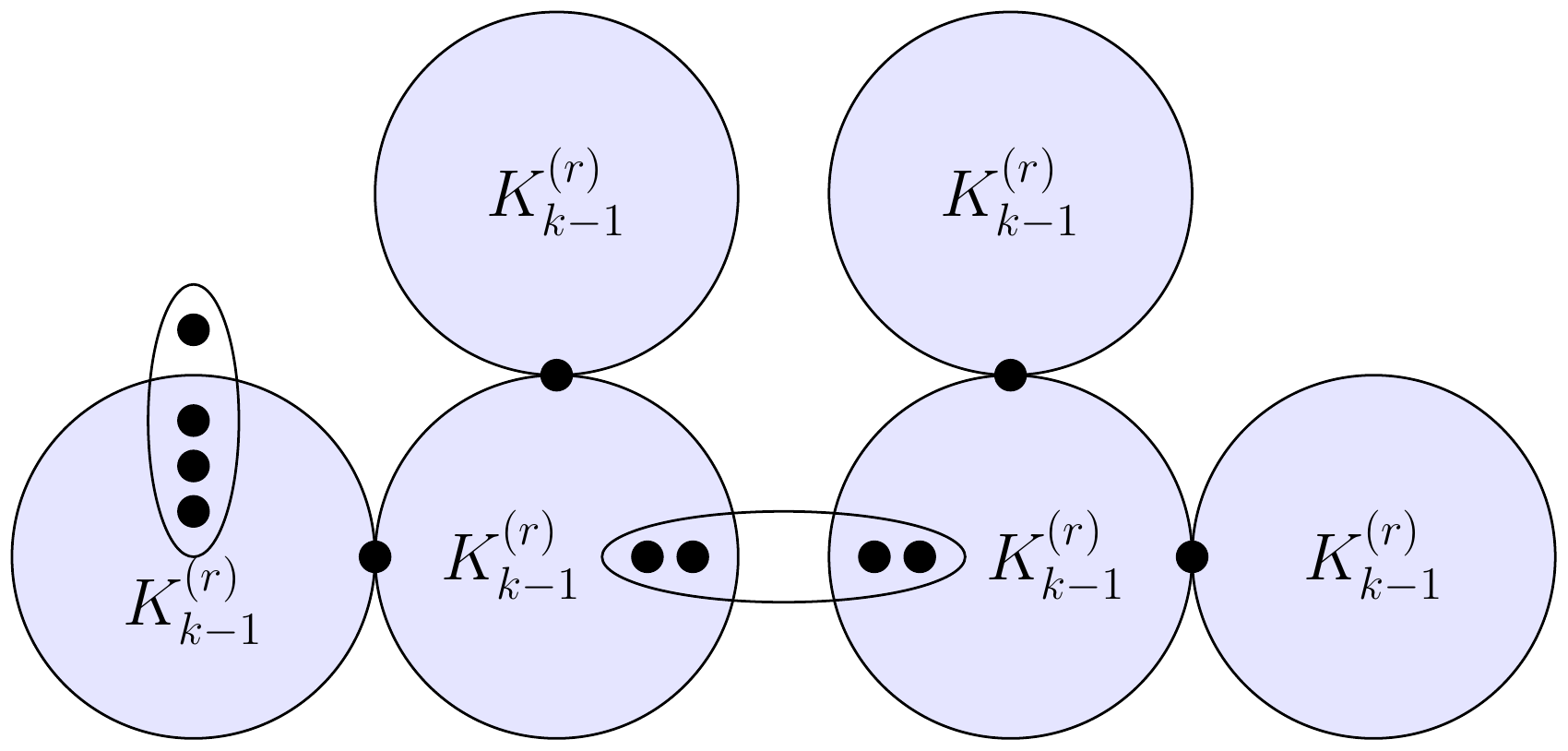}
\caption{An example of a hypergraph from Construction~\ref{const42}.}
\end{figure}
We will see in Section~\ref{sec:proof} that in the case of $m\leq r$  (and $k\geq r+4\geq 7$) these are the only extremal hypergraphs.

\section{\bf  Representative pairs, the structure of Berge $F$-free hypergraphs}\label{sec:tools}

In this section we collect some tools and statements  developed and used in~\cite{FKL}. We do not repeat their proofs.

\begin{definition}\label{def51}
For a hypergraph $\mathcal H$, a {\bf system of distinct representative pairs (SDRP) of $\mathcal H$} is a set of distinct pairs $A = \{\{x_1, y_1\}, \ldots, \{x_s, y_s\}\}$ and a set of distinct hyperedges $\cA = \{f_1, \ldots f_s\}$ of $\mathcal H$ such that for all $1\leq i \leq s$

${}$\quad ---  $\{x_i, y_i\} \subseteq f_i$, and\\
${}$\quad ---  $\{x_i, y_i\}$ is not contained in any $f \in \mathcal H - \{f_1, \ldots, f_s\}$.
\end{definition}

 \begin{lem}\label{le:maxA} Let $\mathcal H$ be a  hypergraph, let $(A, \cA)$ be an \SDRP of $\mathcal H$ of maximum size.
Let $\cB:= \cH \setminus \cA$ and let ${B} = \partial_2 \cB$ be the $2$-shadow of ${\mathcal B}$.
For a subset $S \subseteq B$, let ${\mathcal{B}_S}$ denote the set of hyperedges that contain at least one edge of $S$. Then for all nonempty $S \subseteq B$, $|S|< |{\mathcal B}_S|$.
 \end{lem}

Note that $|\mathcal H| = |A| + |{\mathcal B}|$.

\begin{lem}\label{le:BergeFinG}
Let $\mathcal H$ be a hypergraph and let $(A, \cA)$ be an \SDRP of $\mathcal H$ of maximum size.
Let $\cB:= \cH \setminus \cA$,
 ${B} = \partial_2 \cB$, and let $G$ be the graph on $V(\cH)$ with edge set $A\cup B$.
If $G$  contains a copy of a graph $F$, then $\mathcal H$ contains a Berge $F$ on the same base vertex set.
\end{lem}

In this paper, we only use the previous lemma in the case that $F$ is a cycle or path. I.e., if the longest Berge cycle (path) in $\cH$ is of length $\ell$, then the longest cycle (path) in $G$ is also of length at most $\ell$.

{\bf Definition}: For a natural number $\alpha$ and a graph $G$, the \emph{$\alpha$-disintegration} of
a graph $G$ is the  process of iteratively removing from $G$ the vertices with degree
at most $\alpha$  until the resulting graph has minimum degree at least $\alpha + 1$ or is empty.
This resulting subgraph $H(G, \alpha)$ will be called the $(\alpha+1)$-{\em core} of $G$. It is well known (and easy)
that $H(G, \alpha)$ is unique and does not depend on the order of vertex deletion.

The following theorem is a consequence of Kopylov's Theorem~\cite{Kopy} on the structure of graphs without long cycles. We state it in the form that we need.

\begin{thm}[\cite{Kopy}, see also Theorem 5.1 in~\cite{FKL}]
\label{le:Kopy}
Let $k \geq 5$ and let $t = \lfloor \frac{k-1}{2}\rfloor$.
Suppose that $G$ is an $n$-vertex graph with no  cycle of length at least $k$.
If $G$ is $2$-connected and $n\geq k$ then there exists a subset
$S\subset V(G)$, $s:=|S|$,
  $k-t\leq s\leq k-2$ (i.e., $2\leq k-s\leq t$),
such that the vertices of $V\setminus S$ can be removed by a
$(k-s)$-disintegration.
\end{thm}

\begin{lem}[\cite{FKL},  Lemma 5.3]\label{lem:shadow}
Let $w, \, r\geq 2$ and let $\cH$ be a $w$-vertex $r$-graph.
Let $\overline{\partial_2 \cH}$ denote the family of pairs of $V(\cH)$ not contained in any member of $\cH$ (i.e., the complement of the $2$-shadow).
Then
\[
     |\cH| +    | \overline{\partial_2 \cH}| \leq 
\left\{\begin{array}{ll}\dbinom{w}{2}        &\mbox{for  $2\leq w\leq r+2$},\\
   \dbinom{w}{r}   &\mbox{for $r+2\leq w $}. \end{array} \right.
\]
Moreover,  equality holds if and only if \\
${}$\quad ---  $w > r+2$ and $\mathcal H$ is complete, or\\
${}$\quad ---  $w = r+2$ and either $\mathcal H$ or $\overline{\partial_2 \cH}$ is complete.
  \end{lem}
  
We say a graph $G$ is {\em hamilton-connected} if for any $x,y \in V(G)$, $G$ contains a path from $x$ to $y$ that covers $V(G)$. 
\begin{lem}[\cite{FKLhamcon},  Theorem 5]\label{hamcon}
Let $G$ be an $n$-vertex graph with minimum degree $\delta(G) \geq 2$. If $e(G) \geq {n-1 \choose 2} + 2$ then $G$ is hamilton-connected unless $G$ is obtained from $K_{n-1}$ by adding a vertex
 of degree 2. 
\end{lem}

\section{Maximal mixed hypergraphs} \label{sec:2r}

One of our tools is the notion of {\em mixed hypergraphs}.
For $r\geq 3$, a {\em $(2,r)$ mixed hypergraph} is a triple $\cM=(A, \cB,V)$, where $V$ is a vertex set, $A$ is the edge set of a graph, 
$\cB$ is an $r$-graph (i.e.,
$A\subseteq {V \choose 2}$, $\cB\subseteq {V \choose r}$) such that $A\cup \cB$  satisfies the {\em Sperner property}: there is no $a\in A$, $b\in \cB$ with $a\subset b$.
We often will denote the 2-shadow $\partial_2 \cB$ by $B$.

Let $m_r(n,k)$ denote the maximum size of a mixed hypergraph $\cM$ on $n$ vertices such that $\partial_2\cM$ does not contain any cycle of length $k$ or longer.

\begin{lem}\label{lem:61}
\[
   \EG_r(n,k)\leq m_r(n,k).
 \]
 \end{lem}

\noindent
{\em Proof.}\enskip
Let $\mathcal H$ be an $r$-uniform hypergraph on $n$ vertices with no Berge cycle of length $k$ or longer
 ($k \geq r+3\geq 6$) with $EG_r(n,k)$ edges.
Let $(A, \cA)$ be an \SDRP of $\mathcal H$ of maximum size.
Let $\cB:= \cH \setminus \cA$,
 ${B} = \partial_2 \cB$.
By definition,  $\cM:=(A, \cB,V)$ is a $(2,r)$ mixed hypergraph ($r\geq 3$) with vertex set $V$.
By Lemma~\ref{le:BergeFinG} the graph $G$ with edge set $A\cup B$ does not contain a cycle of length $k$ or longer.
Hence
\[    \EG_r(n,k)=|\cH|=|A|+|\cB|\leq m_r(n,k).  \quad\quad\Box
\]
We will show that these two functions are very close to each other and determine $m_r(n,k)$ for all $n$ (when $k\geq r+4$, $r\geq 3$).
We need more definitions and constructions.

A sequence of sets  $\cS=(V_1, \dots, V_p)$, $V_i\subseteq V$, is called a (linear) {\em hypergraph forest} with vertex set $V$ if
\begin{equation}\label{eq:62}
    | (V_1 \cup \dots \cup V_{i-1}) \cap V_{i}|\leq 1
\end{equation}
holds for each $2\leq i\leq p$.
To avoid trivialities we usually suppose that $|V_i|\geq 2$ for each $i$.
If $\sum_i (|V_i|-1)=|V|-1$ then equality holds in~\eqref{eq:62} for all $i$, and we call $\cS$ a {\em hypergraph tree}.

\begin{const}\label{const:mixed63}
Write $n$ in the form of $(k-2)\lfloor \frac{n-1}{k-2}\rfloor +m$ where $1\leq m\leq k-2$.
Let $p:=\lfloor \frac{n-1}{k-2}\rfloor$.
In case of $m= 1$, let $V_1, \dots, V_p$  be a sequence of $(k-1)$-element subsets of $[n]$ forming a hypergraph tree.
In case of $2\leq m\leq k-2$, let $V_1, \dots, V_{p+1}$ be a sequence of subsets of $[n]$ satisfying~\eqref{eq:62}
 such that  one $V_i$ has $m$ elements and  each other  has $(k-1)$-elements. 
Finally, 
 put either a copy of $K^{(r)}_{|V_i|}$ or $K_{|V_i|}$ 
into each $V_i$.
    \end{const}

Each cycle in the $2$-shadow of any $(2,r)$ mixed family in Construction~\ref{const:mixed63} must be contained in one of the $V_i$'s, so its length is at most $k-1$.
Taking the largest possible mixed hypergraph of this type we get

\begin{equation}\label{eq64}
  m_r(n,k) \geq  f^+_r(n,k):=  \left\lfloor \frac{n-1}{k-2} \right\rfloor \binom{k-1}{r} +
\left\{\begin{array}{ll}\dbinom{m}{2}         &\mbox{for  $1\leq m\leq r+1$},\\
   \dbinom{m}{r}   &\mbox{for $r+2\leq m \leq k-2$}. \end{array} \right.
\end{equation}

\begin{thm}\label{thm:mixed}
Let $r \geq 3$ and $k \geq r+4$, and suppose $\cM$ is an  $n$-vertex $(2,r)$ mixed hypergraph with no cycle of length $k$ or longer in $\partial_2\cM$.
Then $|\cM| \leq f^+_r(n,k) $.
Moreover, equality is achieved if and only if $\cM$ has the structure described in Construction~{\rm \ref{const:mixed63}} above.
\end{thm}

\noindent{\bf Remark.}\quad
This is one point that does not hold for $k=r+3$, because in that case every $\SDRP$ is simply a graph, $\cB=\emptyset$, and according to Kopylov's 
Theorem~\ref{le:Kopy}, there are more extremal graphs than in Construction~\ref{const:mixed63}.

\section{Inequalities} \label{sec:inequalities}

Let $k \geq 5$ and let $t = \lfloor \frac{k-1}{2}\rfloor$, $r\geq 3$, and $k\geq r+3$.
In this section most of the time we suppose that $k\geq r+4$, but almost all inequalities hold for the case $k=r+3$, too. 

Let $\cM=(A, \cB,V)$ be a $(2,r)$ mixed hypergraph such that  $G:=A\cup B$ is an $n$-vertex graph with no  cycle of length at least $k$.

Suppose that $A\cup B$ is $2$-connected and $n\geq k$.
Theorem~\ref{le:Kopy} implies that for some $k-t\leq s\leq k-2$ (i.e., $2\leq k-s\leq t$) there exist an $s$-element set $S\subset V$ such that
\begin{equation}\label{eq101}
\text{\em the vertices of\,\, $A\cup B\setminus S$\,\,  can be removed by a
$(k-s)$-disintegration. 
}
\end{equation}

For the edges of $A$ and $\cB$ contained in $S$ we use Lemma~\ref{lem:shadow} to see that
  \[    |A[S]|+ |\cB[S]|\leq \max\left\{\binom{s}{2}, \binom{s}{r} \right\}.
\]
In the $(k-s)$-disintegration steps, we iteratively remove vertices with degree at most $(k-s)$ until we arrive to $S$.
When we remove a vertex $v$ with degree $\ell \leq (k-s)$ from $G$, $a$ of its incident edges are from $A$, and the remaining $\ell-a$ incident edges eliminate at most $\binom{\ell-a}{r-1}$ hyperedges from $\cB$ containing $v$. Therefore $v$ contributes at most $a + {\ell-a \choose r-1}$ to $|A|+|\cB|$.
Since the function $a+ \binom{\ell-a}{r-1}$ is convex (for nonnegative integers $a$) it takes its maximum at either $a=0$ or $a=\ell$, and since $\ell\leq k-s$ we obtain that
\begin{equation}\label{eq102}
    |A|+|\cB|\leq u_r(n,k,s):=\max\left\{ \binom{s}{2}, \binom{s}{r} \right\} + (n-s) \max\left\{ k-s, \binom{k-s}{r-1}\right\}.
\end{equation}

In the rest of this section we give upper bounds for $u_r(n,k,s)$.
The following inequalities can be obtained by some elementary estimates on binomial coefficients.
The main result of this section is  the following lemma. 

\begin{lem}\label{le103}
If  $r\geq 3$, $k\geq r+4$, $k-t\leq s\leq k-2$, and $n\geq k$, then $u_r(n,k,s)\leq  f_r(n,k)-\dbinom{r}{2}$.
  \end{lem}

{\bf Proof.}
When $s$ is a variable taking only nonnegative integer values, and $r$, $n$ and $k$ are fixed,  the functions $\binom{s}{2}$, $\binom{s}{r}$, $(n-s)(k-s)$ and $(n-s)\binom{k-s}{r-1}$ are convex.
So their maximums and sums, in particular, $u_r(n,k,s)$, are convex, too.
We obtain that
\[    
\max_{k-t\leq s\leq k-2}  u_r(n,k,s) = \max\left\{ u_r(n,k,k-2), u_r(n,k,k-t)   \right\}.
\]

Our first observation is that (for $r\geq 3$, $k\geq r+3$) if $n\geq k-2 +s$, then
\begin{equation}\label{eq104}
         u_r(n,k,s)= u_r(n-k+2, k,s) + (k-2) \max\left\{ k-s, \binom{k-s}{r-1}\right\}.
   \end{equation}


\begin{claim}\label{le104} For $2\leq k-s\leq  t$ and $k\geq r+4$,
\begin{equation}\label{eq105}
 (k-2) \max\left\{ k-s, \binom{k-s}{r-1}\right\} < \dbinom{k-1}{r} -\dbinom{r}{2}.
   \end{equation}
\end{claim}

{\em Proof of Claim~\ref{le104}.}
Because of the convexity of the left-hand side (in variable $s$), it is enough to check the cases $s\in \{ k-t, k-2\}$ (i.e., $k-s\in \{ t, 2\}$, respectively).
We have three cases to consider: when $k-s=2$, when $k-s=t$ and $r\geq k-s$, and finally when $k-s=t$ and $3\leq r\leq t-1$.
Substituting $k-s = 2$ and $k-s = t$ into the left-hand side of~\eqref{eq105}, we get $(k-2)2$ and $(k-2)t$, respectively.
Then (for $k\geq 7$) we have
\[   (k-2)t < \dbinom{k-1}{3}-\dbinom{k-4}{2}\leq   \dbinom{k-1}{r}-\dbinom{r}{2}.
\]
This settles the first two cases.

In the case $k-s=t$ and $r< k-s$,  we need the following inequality (for $3\leq r< t$):
\[   (k-2)\dbinom{t}{r-1}<   \dbinom{k-1}{r}-\dbinom{r}{2}.
\]
We prove the following stronger inequality  (for $3\leq r< t$), because we will 
   use it again.
\begin{equation}\label{eq106}
   \dbinom{k-t}{r}+ (k-3)\dbinom{t}{r-1}<   \dbinom{k-1}{r}-\dbinom{r}{2}.
   \end{equation}
Since $ \binom{k-t}{r}\geq  \binom{t+1}{r}\geq \binom{t}{r-1}$ (for $t\geq r$), equation~\eqref{eq106} completes the proof of~\eqref{eq105}.

Returning to the proof of~\eqref{eq106}
note that  (since $2\leq r-1\leq t-2$)
\[ \binom{r}{2} < \binom{r+1}{2}\leq \binom{t}{2} \leq \binom{t}{r-1}.
\]
So~\eqref{eq106} is implied by the inequality below.
\begin{equation}\label{eq1060}
   \dbinom{k-t}{r}+ (k-2)\dbinom{t}{r-1}\leq   \dbinom{k-1}{r}.
   \end{equation}
We give a purely combinatorial proof of~\eqref{eq1060}. 

Define four $r$-graphs with vertex set $[k-1]$.
$\cF_0:= \dbinom{[k-1]}{r}$, \enskip $\cF_1:=\dbinom{[k-t]}{r}$, \\
$\cF_2:=\left\{ e\cup \{i\}: e\in \dbinom{[t]}{r-1}, i\in [k-t+1,k-1] \right\}$, and\\
$\cF_3:=\left\{ f\cup \{j\}: f\in \dbinom{[k-t, k-1]}{r-1}, j\in [k-t-1] \right\}$.

Their sizes are $\binom{k-1}{r}$, $\binom{k-t}{r}$, $(t-1)\binom{t}{r-1}$, and $(k-t-1)\binom{t}{r-1}$ respectively. 
We claim that $\cF_1$, $\cF_2$, and $\cF_3$ are disjoint. 
Indeed, $|A\cap [k-t]|\leq r-1$ holds for every $A\in \cF_2\cup \cF_3$, so $A\notin \cF_1$.
Also, if $A\in \cF_2$ then $|A\cap [k-t-1]|=r-1>1$ so $A\notin \cF_3$. 
Since the families $\cF_1$, $\cF_2$, and $\cF_3$ are disjoint subfamilies of $\cF_0$, we have
   $|\cF_1|+|\cF_2|+|\cF_3|\leq |\cF_0|$. This completes the proof of~\eqref{eq1060}. \qed

\begin{claim}\label{cl107}
For $k\leq n\leq 2k-3$ and $k\geq r+4$, one has
$ u_r(n,k,k-2)  \leq f_r(n,k) -\dbinom{r}{2}$.  
\end{claim}

{\em Proof.}\enskip
  We have
\begin{eqnarray*} u_r(n,k,k-2) &=&\max\left\{ \binom{k-2}{2}, \binom{k-2}{r} \right\} + (n-k+2) \max\left\{ 2, \binom{2}{r-1}\right\}\\
     &= &\dbinom{k-2}{r}+ 2(n-k+2).
\end{eqnarray*}
We will show
\[
  \dbinom{k-2}{r}+ 2(n-k+2)\leq \dbinom{k-1}{r} +(n-k+1) -\dbinom{r}{2} \quad  \left(\leq f_r(n,k) -\dbinom{r}{2} \right).
   \]
Since $k\geq 7$, we have
\begin{eqnarray*}
\left( \dbinom{k-2}{r}+ 2(n-k+2)\right) -\left( (n-k+1) -\dbinom{r}{2}\right)
    =\dbinom{k-2}{r}+ (n-k+3)+ \dbinom{r}{2}
    \\
    \leq \dbinom{k-2}{r}+ k+ \dbinom{k-4}{2}    \leq \dbinom{k-2}{r}+ \dbinom{k-2}{2}\leq \dbinom{k-1}{r}.
\end{eqnarray*}

\begin{claim}\label{cl108}
For $k\leq n\leq 2k-t-3$, $r\geq t$ and $k\geq r+4$,
\[
  u_r(n,k,k-t) < \dbinom{k-1}{r} +(n-k+1) -\dbinom{r}{2}.
   \]
\end{claim}
Note that the right-hand side is at most
 $f_r(n,k) -\dbinom{r}{2}$.   

{\em Proof.}\enskip
  We have
\begin{eqnarray*} u_r(n,k,k-t) &=&\max\left\{ \binom{k-t}{2}, \binom{k-t}{r} \right\} + (n-k+t) \max\left\{ t, \binom{t}{r-1}\right\}\\
     &= &\dbinom{k-2}{2}+ (n-k+t)t.
\end{eqnarray*}
Moreover
\begin{multline*}
\left( \dbinom{k-t}{2}+ (n-k+t)t\right) -\left( (n-k+1) -\dbinom{r}{2}\right)
    =\dbinom{k-t}{2}+ (t-1)(n-k+t)+(t-2)+ \dbinom{r}{2}\\
    \leq \dbinom{k-t}{2}+  (t-1)(k-3)+(t-2) + \dbinom{k-4}{2} < \dbinom{k-1}{3}\leq \dbinom{k-1}{r}. \quad \Box
\end{multline*}

\begin{claim}\label{cl109}
For $k\leq n\leq 2k-t-3$, $r< t$ and $k\geq r+4$,
\[
  u_r(n,k,k-t) < \dbinom{k-1}{r}  -\dbinom{r}{2}.
   \]
\end{claim}
Note that the right-hand side is at most
$f_r(n,k) -\dbinom{r}{2}$.  

{\em Proof.}\enskip
  We have
\begin{eqnarray*} u_r(n,k,k-t) &=&\max\left\{ \binom{k-t}{2}, \binom{k-t}{r} \right\} + (n-k+t) \max\left\{ t, \binom{t}{r-1}\right\}\\
     &= &\dbinom{k-t}{r}+ (n-k+t)\binom{t}{r-1}  \leq\dbinom{k-t}{r}+ (k-3)\binom{t}{r-1} .
\end{eqnarray*}
Here the right-hand side is less than $\dbinom{k-1}{r}  -\dbinom{r}{2}$ by \eqref{eq106}.
\qed

\medskip
{\em Proof of Lemma~\ref{le103}.} \enskip
Because of the convexity of $u_r(n,k,s)$  (in the variable $s$), it is enough to check the cases  $k-s\in \{ 2, t\}$.
Note that for $n_1, n_2\geq 2$
\begin{eqnarray}\label{eq1012}
               f_r(n_1,k)+ f_r(n_2,k)&\leq&  f_r(n_1+n_2-1,k)\\
     f_r^+(n_1,k)+ f_r^+(n_2,k)&\leq&  f_r^+(n_1+n_2-1,k)   \label{eq1013}
   \end{eqnarray}
and here equalities hold for $n_2=k-1$.
(If we define $f_r(1,k)=f_r^+(1,k)=0$, then we can use~\eqref{eq1012},~\eqref{eq1013} for these values, too).

If $s=2$ and $k\leq n\leq 2k-3$, then  Claim~\ref{cl107} yields $u_r(n,k,k-2) \leq  f_r(n,k) -\dbinom{r}{2}$.
For $n\geq 2k-2$  we use \eqref{eq104}, then the induction hypothesis $u_r(n-k+2,k, k-2)\leq f_r(n-k+2,k)$, and then Claim~\ref{le104} (equation~\eqref{eq105})  implies that
\begin{eqnarray*}
           u_r(n,k,k-2)&=& u_r(n-k+2, k,k-2) + (k-2) 2
       <  f_r(n-k+2,k)  + \dbinom{k-1}{r} -\dbinom{r}{2}
   \\ &=&  f_r(n-k+2,k)  + f_r(k-1,k)-\dbinom{r}{2}
        =  f_r(n,k) -\dbinom{r}{2},
   \end{eqnarray*}
and we are done.

When $k-s=t$ the proof is similar.
For $k\leq n\leq 2k-t-3$,   Claim~\ref{cl108} and Claim~\ref{cl109} yield $u_r(n,k,k-t)< f_r(n,k) -\dbinom{r}{2}$.
For $n\geq 2k-t-2$,  we use \eqref{eq104}, then the induction hypothesis $u_r(n-k+t,k, k-t)\leq f_r(n-k+t,k)$,  and then Claim~\ref{le104} (equation~\eqref{eq105})
 implies that
\begin{eqnarray*}
                   u_r(n,k,k-t)&=& u_r(n-k+2, k,k-2) + (k-2) \max\left\{ t, \binom{t}{r-1}\right\}\\
       &<&  f_r(n-k+2,k) + \dbinom{k-1}{r}  -\dbinom{r}{2}
   \\ &=&  f_r(n-k+2,k) + f_r(k-1,k) -\dbinom{r}{2}
        =  f_r(n,k) -\dbinom{r}{2}.
   \end{eqnarray*}
 \qed

\section{Proofs of the main results} \label{sec:proof}

In this section we first prove Theorem~\ref{thm:mixed} and then  Theorem~\ref{main:all} for all $n\geq k$ (and $r\geq 3$, $k\geq r+4$).

\subsection{Proof of Theorem~\ref{thm:mixed} about mixed hypergraphs}
Let $\cM=(A, \cB,V)$ be a $(2,r)$ mixed hypergraph such that  $G:=A\cup B$ is an $n$-vertex graph with no  cycle of length at least $k$ ($B:= \partial_2 \cB$ and $A\cap B=\emptyset$).
Let $V_1, V_2, \dots, V_q$ be the vertex sets of the standard (and unique) decomposition of $G$ into blocks of sizes $n_1, n_2, \dots, n_q$.
Then the graph $A\cup B$ restricted to $V_i$, denoted by $G_i$, is either a $2$-connected graph or a single edge (in the latter case $n_i=2$),
 each edge from $A\cup B$ is contained in a single $G_i$, and $\sum_{i=1}^q (n_i-1)\leq (n-1)$.
This decomposition yields a decomposition of $A=A_1\cup A_2\cup \dots\cup A_q$ and
  $B=B_1\cup B_2\cup \dots\cup B_q$, $A_i\cup B_i=E(G_i)$.
If an edge $e\in B_i$ is contained in $f\in \cB$, then $f\subseteq V_i$ (because $f$ induces a
$2$-connected graph $K_r$ in $B$), so the block-decomposition of $G$ naturally extends to $\cB$,
  $\cB_i:= \{ f\in \cB: f\subseteq V_i\}$ and we have $\cB= \cB_1\cup \dots \cup \cB_q$, and $B_i=\partial_2 \cB_i$.
By definition, $G$ has no cycle of length $k$ or longer, so the same is true for each $G_i$.
Suppose that the size of $A\cup \cB$ is as large as possible, $\cM$ is extremal, $|\cM|=m_r(n,k)$.

Lemma~\ref{lem:shadow} implies that for $n_i\leq k-1$,
  \[|A_i|+|\cB_i|\leq \max\left\{ \dbinom{n_i}{2} ,\dbinom{n_i}{r}\right\}= f^+_r(n_i, k),\]
and equality holds only if $A_i$ is the complete graph (and $\cB_i=\emptyset$) or $\cB_i$ is the $r$-uniform complete $r$-graph (and $A_i=\emptyset$).

Lemma~\ref{le103} implies that in the case $n_i\geq k$
\begin{equation}\label{eq111}
   |A_i|+|\cB_i|\leq  f_r(n_i,k)-\dbinom{r}{2} < f^+_r(n_i, k).
   \end{equation}
Adding up these inequalities for all $1\leq i\leq q$ and applying~\eqref{eq1013}, we get
 \begin{equation}\label{eq112}
     \sum_i (|A_i|+|\cB_i|)\leq  \sum_i  f^+(n_i,k) \leq f_r^+(1 + \sum_i (n_i-1),k)\leq f_r^+(n,k).
\end{equation}
Since $ f_r^+(n,k)\leq m_r(n,k)$, here equality holds in each term.
Consequently $n_i< k$ for each $i$, and all but at most one of them should be $k-1$.
Otherwise we can use the inequality
\[ \max\left\{ \dbinom{a}{2}, \dbinom{a}{r}\right\} +\max\left\{ \dbinom{b}{2}, \dbinom{b}{r}\right\} <
      \max\left\{ \dbinom{a-1}{2}, \dbinom{a-1}{r}\right\} +\max\left\{ \dbinom{b+1}{2}, \dbinom{b+1}{r}\right\}
\]
which holds for all $1< a\leq b< k-1$ (and $3\leq r$, $r+4\leq k$).
(The inequality  $f(a)+ f(b)\leq f(a-1)+f(b+1)$ holds for every convex function $f$, and here equality holds only if the four points
$(a-1,f(a-1))$, $(a,f(a))$, $(b,f(b))$, and $(b+1,f(b+1))$ are lying on a line).
So $\cM$ is a linear tree formed by cliques, as described in Construction~\ref{const:mixed63}.
\qed

\subsection{Proof of Theorem~\ref{main:all} for $m> r+1$}
Let $\mathcal H$ be an $r$-uniform hypergraph on $n$ vertices with no Berge cycle of length $k$ or longer ($r\geq 3$,  $k \geq r+4$).
Suppose that $|\cH|$ is maximal, $|\cH|= \EG_r(n,k)$. We have $f_r(n,k)\leq  \EG_r(n,k)$ by Constructions~\ref{const41} and~\ref{const42}.

Let $(A, \cA)$ be an \SDRP of $\mathcal H$ of maximum size.
Let $\cB:= \cH \setminus \cA$, ${B} = \partial_2 \cB$.
By Lemma~\ref{le:BergeFinG} the graph $G$ with edge set $A\cup B$ does not contain a cycle of length $k$ or longer.
In other words, $\cM=(A, \cB,V)$ is a $(2,r)$ mixed hypergraph such that  $G:=A\cup B$ is an $n$-vertex graph with no  cycle of length at least $k$.
Then Theorem~\ref{thm:mixed} implies that
    \begin{equation}\label{113}   |A|+|\cB|\leq f^+_r(n,k).
      \end{equation}
Since $n=(k-2)p +m$ where $1\leq m\leq k-2$ and $m\geq r+2$ we have $f^+_r(n,k)=f_r(n,k)$
by~\eqref{eq:32} and~\eqref{eq64}. We obtained that  $\EG_r(n,k)=f_r(n,k)$, as claimed.

Equality can hold in~\eqref{113} only if $\cM$ has the clique-tree structure, $V_1, V_2, \dots, V_{p+1}$, described in Construction~\ref{const:mixed63}.
In the case of $m\geq r+3$ each block is a complete $r$-uniform hypergraph, so Construction~\ref{const:mixed63} and Construction~\ref{const41} coincide, and we are done.

In the case $m=r+2$, Theorem~\ref{thm:mixed} implies that all but one block define complete $r$-graphs and for one of them, say $V_\ell$, $\cM|V_\ell$ could be either $K_{r+2}$ or $K_{r+2}^{(r)}$. 
If $\cM|V_\ell=K_{r+2}^{(r)}$, then  $\cH|V_\ell= K_{r+2}^{(r)}$,  so $\cA=\emptyset$, $\cB=\cH$ and we are done.
Consider the other case, $\cM|V_\ell= K_{r+2}^{(2)}$, i,e., $A= G|V_\ell$ is a complete graph (and $\cB= \cup_{i\neq \ell} K_{k-1}^{(r)}[V_i]$).
We claim that $\cH|V_\ell=K_{r+2}^{(r)}$ which completes the proof in this subsection.

Suppose, on the contrary, that there exists an $f_i\in \cA$ such that $\{x_i,y_i\}\subset V_\ell$, $\{x_i, y_i, z_i\}\subset f_i$
 such that $z_i\notin V_\ell$.
One of the pairs of $x_iz_i$ and $y_iz_i$ is not an edge of $G$, say it is $x_iz_i$.
Then removing $x_iy_i$ from $A$ and replacing it by $x_iz_i$, one obtains an \SDRP $A'$ ($\cA$ and $\cB$ are unchanged).
In this case, $E(G')= E(G)\setminus \{ x_iy_i\} \cup \{ x_iz_i\} $ has a different structure (not a tree of cliques), so it could not be optimal by Theorem~\ref{thm:mixed}.
Therefore such $f_i$ does not exist, i.e., $f_i\subset V_\ell$. In other words $\cA\subseteq K_{k-1}^{(r)}[V_i]$.
Since $|A|= \binom{r+2}{2}=\binom{r+2}{r}$,  $\cA$ is a complete $r$-graph on  $V_\ell$. \qed

\subsection{Proof of Theorem~\ref{main:all} for $m\leq r+1$, preparations}

This is a continuation of the previous two subsections.

Consider an extremal $\cH$ (i.e., $|\cH|=\EG_r(n,k)\geq f_r(n,k)$) and the blocks $G_1, \dots, G_q$ of $G$ with vertex sets $V_1, \dots, V_q$ where $|V_i|=n_i\geq 2$.
As we have seen in~\eqref{eq111} and~\eqref{eq112}, 
\begin{equation}\label{eq1140}
    |\cH|=\sum_i \left(   |A_i|+|\cB_i| \right)\leq f^+_r(n,k) -\binom{r}{2}
  \end{equation}
if for any $i$, $n_i \geq k$.
For $m=r+1$, here the right-hand side is
\[ p\binom{k-1}{r}+ \binom{m}{2}-\binom{r}{2} < p\binom{k-1}{r} +\binom{r+1}{r}= f_r(n,k).
   \]
Similarly in the case $m\leq r$, the right-hand side is
\[ p\binom{k-1}{r}+ \binom{m}{2}-\binom{r}{2} < p\binom{k-1}{r} +m-1= f_r(n,k).
   \]
So from now on, we may suppose that $n_i\leq k-1$ for all $i$.

\begin{claim}\label{114} There are exactly $p$ blocks $V_i$ of size $k-1$, $n_i=k-1$.
    \end{claim}

{\em Proof.}\enskip
For $1\leq x\leq k-1$, define $f(x):= f_r^+(x,k)=\max\{ \binom{x}{2}, \binom{x}{r}\}$.
Let $f(x_1, \dots , x_q):= \sum_i f(x_i)$.
We want to estimate $f(n_1, \dots, n_q)$, so define $x_i:=n_i$.
Let $n':= 1+ \sum_i(n_i-1)$; we have $n'\leq n$.
In case of $2\leq x_i\leq x_j< k-1$ we are going to replace $x_i$ by $x_i-1$ and $x_j$ by $x_j+1$.
During this process $f$ never decreases and it ends when all but one  $x_i$'s become $1$ or $k-1$.
Then the value of $f$ is exactly $f_r^+(n', k)$ and since $\sum_{1\leq i\leq q}  x_i= \sum_i n_i = n'+q-1$ is unchanged,
  in the last step our sequence contains $(k-1)$ exactly $p$ times.

If the number of $(k-1)$'s is unchanged, then there is nothing to prove.
Otherwise, after some step the pair $x$ and $k-2$ ($2\leq x\leq k-2$) is replaced by $(x-1)$ and $(k-1)$.
Then the value of $f$ increased by $f(k-1)+f(x-1)-f(k-2)-f(x)$.
Since $f(k-1)=\binom{k-1}{r}$ and $f(k-2)=\binom{k-2}{r}$ the increment is
\[  \binom{k-1}{r}+ \max\{ \binom{x-1}{2}, \binom{x-1}{r}\}-  \binom{k-2}{r}-\max\{ \binom{x}{2}, \binom{x}{r}\}.
\]
This is at least
\[
\binom{k-2}{r-1}- \max\{ x-1, \binom{x-1}{r-1}\}\geq \binom{k-2}{r-1}-\binom{k-3}{r-1}=\binom{k-3}{r-2}\geq \binom{r+1}{r-2} >\binom{r}{2}. 
\]
In this case $|\cH|< f_r^+(n', k)-\binom{r}{2}\leq f_r(n',k)\leq f_r(n,k)$, a contradiction.
\qed

\begin{claim}\label{115} If a block $V_i$ is of size $k-1$, then $e(B_i)\geq \binom{k-2}{2}+r-1$.
    \end{claim}

{\em Proof.}\enskip
If $|\cB_i|> \binom{k-2}{r}$ then the Kruskal-Katona Theorem (or a simple double counting) implies that $|\partial_2 \cB_i|\geq \binom{k-2}{2}+r-1$, and we are done.

If  $|\cB_i|\leq  \binom{k-2}{r}$ then we  use Lov\'asz' version of the Kruskal-Katona theorem.
Write $|\cB_i|$ in the form of $\binom{x}{r}$, where $x$ is a real number $0\leq x \leq k-2$ and  (only in this paragraph) $\binom{x}{r}$ is defined as the {\em real} polynomial $x(x-1)\dots (x-r+1)/ r!$ for $x\geq r-1$ and $0$ otherwise.
We obtain $|\partial_2\cB_i|\geq \binom{x}{2}$. Since $A_i$ and $B_i$ are disjoint, we have $|A_i|\leq \binom{k-1}{2}-\binom{x}{2}$. 
So,
\begin{equation}\label{eq11401}
 |A_i|+ |\cB_i|\leq \binom{k-1}{2}-\binom{x}{2}+ \binom{x}{r}
\end{equation}
holds for some $0\leq x \leq k-2$. In this range the right-hand side (as a polynomial of variable $x$) is maximized at $x=k-2$. 
Hence~(\ref{eq11401}) yields
\[  |A_i|+ |\cB_i|\leq \binom{k-1}{2}-\binom{k-2}{2}+ \binom{k-2}{r}.
\]
Here the right-hand side is less than $\binom{k-1}{r}-\binom{r}{2}$ which (as we have seen in~\eqref{eq1140}) leads to the contradiction
 $|\cH|< f_r(n,k)$.
\qed

\begin{claim}\label{116} If a block $V_i$ is of size $k-1$, then $\cB_i= K_{k-1}^{(r)}$, a complete $r$-graph.
    \end{claim}

{\em Proof.}\enskip
Suppose that there exists an $r$-set $f\subset V_i$, $f\notin \cH$.
Consider the hypergraph $\cH\cup \{ f\}$.
By the maximality of $\cH$, $\cH \cup \{f\}$  contains a Berge cycle $C$ of length at least $k$, say with base vertices $\{v_1, \ldots, v_\ell\}$ and edges $\{f_1, \ldots, f_\ell\}$ where $f_\ell = f$ (and so $v_1, v_\ell \in V_i$). Since $|V_i|= k-1$, there is a base vertex of $C$ not contained in $V_i$. Therefore we may pick a segment $P$ of $C$ (a Berge path in $\cH$) say $\{v_{a}, v_{a +1}, \ldots, v_{b}\}$, $\{f_a, \ldots, f_{b-1}\}$ such that $v_a, v_b \in V_i$ but $\{v_{a+1}, \ldots, v_{b-1}\} \cap V_i = \emptyset$.

Since each $r$-edge in $\cB_i$ yields a clique of order $r$ in $B_i$, we have  $\delta(B_i) \geq r-1 \geq 2$. By Claim~\ref{115} and Lemma~\ref{hamcon}, $B_i$ is hamilton-connected unless $r = 3$ and $B_i$ is a clique on $k-2$ vertices with a vertex $x$ of degree 2. If the latter holds, then for a neighbor $y$ of $x$, the edge $xy$ is contained in exactly one triangle in $B_i$. But then $xy$ can only be contained in one $r$-edge of $\cB$, contradicting Lemma~\ref{le:maxA}. So we may assume $B_i$ has a hamilton path between any two vertices, in particular by Lemma~\ref{le:BergeFinG}, there is a Berge path $P'$ of length $k-2$
from $x_a$ to $x_b$ containing all $k-1$ vertices of $V_i$ as base vertices and using only the edges from $\cB_i$. 
The cycle $P\cup P'$ is a Berge cycle in $\cH$ of length at least $k$, a contradiction.
Therefore such an edge $f$ cannot exist, $\cH|V_i= K_{k-1}^{(r)}$.

Finally, there is no $A$-edge in $V_i$. If $\{ x,y\}\subset V_i$ is an $A$-edge, then no $\cB$-edge can contain $\{x,y\}$. So all the $\binom{k-3}{r-2}\, (\geq k-3)$  subsets of $V_i$ of size $r$ and containing $xy$ should belong to $\cA$.
Therefore $V_i$ must contain at least as many $A$-edges. But $|A_i|\leq k-1-r (={k-1 \choose 2} - {k-2 \choose 2} -r + 1)$ by Claim~\ref{115}.
\qed

\subsection{Proof of Theorem~\ref{main:all} for $m\leq r+1$, the end}

This is a continuation of the previous three subsections.

Consider an extremal $r$-graph $\cH$ on the $n$-element vertex set $V$  (i.e., $|\cH|=\EG_r(n,k)\geq f_r(n,k)$) where $n=p(k-2)+m$, $1\leq m\leq r+1$.
Using the previous subsection, we define a different split of $\cH$.

Let $\cV:= \{ V_i: |V_i|=k-1\}$. By Claim~\ref{114}, $|\cV|=p$.
Let $H$ be the graph whose edge set is the union of the complete graphs on $V_i\in \cV$, so $|E(H)|=p\binom{k-1}{2}$ and it has a forest like structure of cliques
 (i.e., every cycle in $H$ is contained in some $V_i\in \cV$).
Let $C_1, \dots, C_m$ be the vertex sets of the connected components of $H$.
The graph $H$ necessarily consists of $m$ (nonempty) components, $\cup C_\alpha=V$ ($1\leq \alpha\leq m$), some of them could be singletons.
Let $H_\alpha:= H|C_\alpha$, $\cH_\alpha:= \cup \{ \cB_i: V_i\in \cV, V_i\subset C_\alpha\}$, and $\cD:= \cH\setminus (\cup \cH_\alpha)$.
Note that every edge of $H$ used to be a $B$-edge, $\cH_\alpha\subseteq \cB$ for all $1\leq \alpha\leq m$, and $\cD$ is the set of edges in $\cH$ not contained in some $K_{k-1}^{(r)}$.

Our main observation is the following which is  implied by Claim~\ref{116}.

\begin{claim}\label{claim117}
If $x,y\in C_\alpha$, $x\neq y$ then there exists  an $x$-$y$ Berge path of length at least $k-2$ consisting only of $\cH_\alpha$ edges.
Moreover, if $xy\notin E(H_\alpha)$ then there exists  such a path of length at least $2k-4$. \qed
   \end{claim}

Suppose that $f,f'\in \cD$, ($f\neq f'$), $x_\alpha\in C_\alpha \cap f$, $x'_\alpha\in C_\alpha \cap f'$, $x_\beta\in C_\beta \cap f$,
  and $x'_\beta\in C_\beta \cap f'$, ($\alpha \neq \beta$), then
\begin{equation}\label{eq118}
     x_\alpha=x'_\alpha \text{ and }x_\beta=x'_\beta.
   \end{equation}
For example, if $x_\alpha\neq x'_\alpha$ and $x_\beta\neq x'_\beta$, then there is a Berge path $P_\alpha$
  of length at least $(k-2)$ connecting $x_\alpha$ with $x_\alpha'$,
 $P_\alpha\subset \cH_\alpha$ and  another Berge path $P_\beta$ of length at least $(k-2)$
 connecting $x_\beta$ with $x_\beta'$, $P_\beta\subset \cH_\beta$, and these,
 together with $f$ and $f'$ form a Berge cycle of length at least $2k-2$, a contradiction.
The case $|\{ x_\alpha, x_\beta\}\cap \{x'_\alpha, x'_\beta\}|=1$ is similar: we find a Berge cycle in $\cH$ of length at least $k$.  \qed

The same proof, and the second half of Claim~\ref{claim117} imply that
\begin{equation}\label{eq119}  \partial_2 \cH| V_\alpha = H_\alpha.
\end{equation}
In other words, if $f\in \cH\setminus \cH_\alpha$ then
\[
\text{$|f\cap C_\alpha|\geq 2$ implies that $\exists V_i\in \cV$, $V_i\subseteq C_\alpha$ such that $C_\alpha\cap f= V_i\cap f$.}
\]
Indeed, otherwise there are $x,y\in f$ and a  Berge $x,y$-path in $\cH_i$ of length at least $2k-4$, which together with $f$ form a
 Berge cycle of length at least $2k-3$.

For a subset $S\subseteq V$, define $\varphi(S)$ as the set of indices $1\leq \alpha \leq m$ for which $S\cap C_\alpha \neq \emptyset$. 
Equation~\eqref{eq118} can be restated as follows 
\begin{equation}\label{eq1110}\text{if $\{ \alpha, \beta\}\subseteq \varphi(f) \cap \varphi(f')$ then $C_\alpha\cap f= C_\alpha\cap f'$ is a singleton,}
\end{equation}
and similarly for $\beta$. This implies that $\varphi(f)\neq \varphi(f')$ for $f\neq f'$, $f,f'\in \cD$.
Even more, the family $\{\varphi(f): f\in \cD \}$ has the Sperner property.
This means that for  $f,f'\in \cD$ with $f\neq f'$, one cannot have $\varphi(f) \subsetneq \varphi(f')$.
Indeed, $|\varphi(f)|< r$ implies that there exists a $C_\alpha$ with $|C_\alpha \cap f|\geq 2$,
 equation~\eqref{eq119} implies that $|\varphi(f)|\geq 2$ for every $f\in \cD$, so there exists a
  $\beta\in \varphi(f)$, $\alpha\neq \beta$.
But then $\{ \alpha, \beta\} \subseteq \varphi(f)\cap \varphi(f')$ and~\eqref{eq1110} implies that $|C_\alpha\cap f|=1$, a contradiction.

The following claim on the intersection structure of the edges in $\cD$ is a generalization of~\eqref{eq1110} which can be considered as the case $\ell=2$.
(Technically, two hyperedges sharing at least two vertices form a Berge cycle of length $2$.)
\begin{claim}\label{claim1111}
 Let $\cF=:\varphi(\cD)=\{ \varphi(f): f\in \cD\}$.
Suppose that $\{ \alpha_1, \dots, \alpha_\ell\} \subset \{ 1, \dots, m\}$ and $\varphi(f_1), \dots, \varphi(f_\ell)\in \cF$
form a Berge cycle in $\cF$.
Then for each $j$,  the sets $C_{\alpha_j}\cap f_j= C_{\alpha_j}\cap f_{j-1}$ are singletons.
  \end{claim}
{\bf Proof.}\enskip
Otherwise, we can relabel $j:=1$ and  find two distinct vertices $x_1$ and $x_1'$ such that $x_1\in C_{\alpha_1}\cap \varphi(f_1)$ and $x_1'\in  C_{\alpha_1}\cap \varphi(f_\ell)$. 
Furthermore, let $x_i, x_i'\in C_{\alpha_i}$ such that $\{ x_{i-1}, x_i\} \subset f_i$ for all $1\leq i\leq \ell$ ($x_0:=x_\ell$, etc.),  $P_i$ a Berge path in $\cH_i$ connecting
 $x_i$ with $x_i'$. These paths could be empty (if $x_i=x_i'$) but  by~Claim~\ref{claim117} we can choose $P_1$ so that its length is at least $k-2$.
Then $f_1,P_1, f_2, P_2, \dots, f_\ell, P_\ell$ form a cycle of length at least $k$, a contradiction.  \qed

\bigskip
{\bf Case 1: }there exists an $f$ such that $|\varphi(f)|=r$. 
Then $m\geq r$, so $m\in \{r,r+1\}$.
If $m=r$, then (because of the Sperner property) $|\cD|=1 < m-1$, a contradiction.
So assume $m=r+1$. Let $\alpha:= [m]\setminus \varphi(f)$. We have $\alpha\in \varphi(f')$ for all other $f'\in \cD$.
Since $|\cD|\geq r+1>3$, there are at least two more $f_2\neq f_3\in \cD\setminus \{ f\}$.

Consider first the case that $|C_\alpha\cap \varphi(f_2)|\geq 2$ for some $f_2\in \cD$.
The Sperner property implies that there are distinct $\alpha_2, \alpha_3\in [m]\setminus \alpha$ such that $\alpha_2\in \varphi (f_2)\setminus \varphi(f_3)$ and
$\alpha_3\in \varphi (f_3)\setminus \varphi(f_2)$. Then $\alpha, \alpha_2, \alpha_3$ with the hyperedges $\varphi(f_2)$, $\varphi(f)$, and $\varphi(f_3)$ form a Berge cycle.
However this cycle does not satisfy Claim~\ref{claim1111}.
So from now on, we may suppose that  $|C_\alpha\cap \varphi(f')|=1$ for all  $f'\in \cD\setminus \{ f\}$.

Suppose that there exists an $f_2\in \cD$ and an $\alpha_2\in [m]$ such that $|C_{\alpha_2}\cap \varphi(f_2)|\geq 2$ (necessarily $\alpha_2\neq \alpha$).
Again  Sperner property implies that there is an $\alpha_3\in [m]\setminus \alpha$ such that
$\alpha_3\in \varphi (f_3)\setminus \varphi(f_2)$ (so we have $\alpha_3\neq \alpha_2$).
Then $\alpha, \alpha_2, \alpha_3$ with the hyperedges $\varphi(f_2)$, $\varphi(f)$, and $\varphi(f_3)$ form a Berge cycle.
However this cycle does not satisfy Claim~\ref{claim1111}.
So from now on, we may suppose that  $|C_{\alpha'}\cap \varphi(f')|=1$ for all  $f'\in \cD$ and all $\alpha'\in [m]$.

Since $|\cD|\geq r+1$ and $[m]$ has exactly $r+1$ $r$-subsets,  $\varphi(\cD)$ is a complete $r$-graph.
Its hyperedges form many Berge cycles, so  Claim~\ref{claim1111} implies that $\cD$ itself is isomorphic to $K_{r+1}^{(r)}$.
Thus $\cH$ is as in Construction~\ref{const41}.

\bigskip
{\bf Case 2:} $|\varphi(f)|<r$ for all $f\in \cD$. 
In this case every $f\in \cD$ has an $\alpha(f)\in [m]$ such that $|C_{\alpha(f)}\cap f|\geq 2$.
For every $f \in \cD$, choose another element $\beta(f)\in \varphi(f)$ ($\beta(f)\neq \alpha(f)$) and consider the graph $T:=\{\{ \alpha(f), \beta(f)\}: f\in \cD \}$.
By Claim~\ref{claim1111} the graph $T$ has no cycle, and the maximality of $|\cH|$ implies that $e(T)=|\cD|\geq m-1$. So $T$ is a tree.
Since $T$ is a tree, 
one cannot replace an edge $\{ \alpha(f), \beta(f)\}$ by the  3-edge $\{ \alpha(f), \beta(f), \gamma(f)\}$ without creating a cycle in the resulting hypergraph and thus violating
Claim~\ref{claim1111}. So $\varphi(f)=  \{ \alpha(f), \beta(f)\}$, and~\eqref{eq119} implies that the structure of $\cD$ is as  in
Construction~\ref{const42}. This completes the proof of Theorem~\ref{main:all}.
\qed

\end{document}